\documentclass{amsart}
\usepackage{amsfonts,amsthm,amsmath}

\newtheorem{axiom}{Axiom}[section]
\newtheorem{defin}{Definition}[section]
\newtheorem{theo}{Theorem}[section]
\newtheorem{propo}{Proposition}[section]
\newtheorem{lemma}{Lemma}[section]
\newtheorem{corol}{Corollary}[section]
\newtheorem{remark}{Remark}[section]

\newcommand{\bax}{\begin{axiom}}
\newcommand{\eax}{\end{axiom}}
\newcommand{\bass}{\begin{assumption}}
\newcommand{\eass}{\end{assumption}}
\newcommand{\bdefi}{\begin{defin}}
\newcommand{\edefi}{\end{defin}}
\newcommand{\bth}{\begin{theo}}
\newcommand{\eth}{\end{theo}}
\newcommand{\bprop}{\begin{propo}}
\newcommand{\eprop}{\end{propo}}
\newcommand{\blem}{\begin{lemma}}
\newcommand{\elem}{\end{lemma}}
\newcommand{\bcor}{\begin{corol}}
\newcommand{\ecor}{\end{corol}}
\newcommand{\brem}{\begin{remark}}
\newcommand{\erem}{\end{remark}}
\newcommand{\bpf}{\begin{proof}}
\newcommand{\epf}{\end{proof}}

\newcommand{\field}[1]{\mathbb{#1}}

\newcommand{\R}{\field{R}}

\newcommand{\ra}{\rightarrow}

\newcommand{\pa}{\partial}
\newcommand{\f}{\frac}

\newcommand{\de}{{\delta}}

\newcommand{\ga}{\gamma}

\newcommand{\cN}{{\mathcal N}}

\newcommand{\cM}{{\mathcal M}}

\begin{document}


\title[]{Odd-symplectic group in first order partial differential equations}

\author{L.Sbano}

\address{
Mathematics Institute, University of Warwick, CV4 7AL Coventry, U.K.\\
{\bf e-mail}: Sbano@maths.warwick.ac.uk,
L.Sbano@mclink.it}

\begin{abstract}
In this paper it is shown that the characteristic vector field associated to
a first order PDE:
\[h(x,\nabla z(x))=0\mbox{ with $x\in\R^n$}\]
has the same form of an infinitesimal  generator of an odd-symplectic transformation 
with contact Hamiltonian $h(x,p)$ on the level set $h=0$. We also study under 
which condition such PDE has a characteristic vector field commuting with 
 a generator of an odd-symplectic action.  
\end{abstract}


\vfill
{\it}

\maketitle


\section*{Introduction}
Very recently R. Cushman \cite{Cushman0}  
 presented the notion of {\it odd-symplectic group}. This notion  was formulated 
in \cite{Gelfand} but differently in \cite{Proctor}; more recently  
in \cite{Cushman,Bates}. Such group can be defined on vector spaces with odd dimension. 
Loosely speaking, the odd-symplectic group 
lies between the standard symplectic group and the contact group. More precisely, it 
 interpolates symplectic groups in different dimensions according to the following exact sequence:
\[\dots\ra Sp(2n+2,\R)\ra Sp(2n+1,\R)\ra Sp(2n,\R)\ra\dots\]
 Although the definition of the  odd-symplectic group is natural, at the moment, there are few 
examples of applications. The aim of this paper is to show that  odd-symplectic transformations 
 rise naturally in the geometric theory of first order PDE and in particular in the study of the characteristics 
 of PDEs of the following form:
\begin{equation}
h(x,\nabla z(x))=0,
\label{h=0}
\end{equation}
where $x\in\R^n$ $h(.,.)$ is a smooth function on $\R^{2n}$.\\
Here is an outline of the paper: in the first section the odd-symplectic group is defined. 
In the following sections we recall the theory of characteristics for first order PDEs and 
 we show that the characteristic vector field associated to (\ref{h=0}) coincides with 
an infinitesimal generator of an odd-symplectic transformation. We show also that  
\[h(x,\nabla z(x))=\f{1}{2}<g\nabla z,\nabla z>+<v,\nabla z>+k\]
 with $g$ constant $n\times n$ symmetric matrix, $k\in \R$ and $v\in\R^n$ has a symmetry 
whose infinitesimal generator is given by an element in the Lie algebra of 
$Sp(2n+1,\R)$.
Finally we show that there is an isomorphism 
between certain class of first order PDE and the Lie algebra of the odd-symplectic group; moreover we show 
 the connection of the odd-symplectic group with Hamilton-Jacobi theory and the Eikonal equation.

\section{Odd-symplectic group and contact geometry}
Let us consider a symplectic vector space $(V,\omega)$, $\dim V=2n+2$. The symplectic group 
is:
\begin{equation}
Sp(V,\omega)\doteq\{A\in Gl(V)| A^*\omega=\omega\}.
\label{simplectic}
\end{equation}
The odd-symplectic group $Sp(V,\omega)_{v_0}$ can be defined as follows:
let $v_0\in V$ be a non-zero vector,  then:
\begin{equation}
Sp(V,\omega)_{v_0}\doteq\{A\in Sp(V,\omega)| Av_0=v_0\}
\label{oddsimplectic}
\end{equation}
Sometimes $Sp(V,\omega)_{v_0}$ is denoted by $Sp(2n+1,\R)$. 
This group has been studied in \cite{Bates,Cushman}. It turns out that 
the odd symplectic group arises naturally in contact geometry. Let us follow \cite{Bates}.
Consider a contact manifold $(\cM,\theta)$ ($\dim\cM =2n+1$)
 where locally, in Darboux coordinates, the contact form $\theta$ is:
\begin{equation}
\theta=dz-y_idx^i
\end{equation}
In \cite{Arnol'd} it is shown that 
one can symplectify $(\cM,\theta)$ into $(\cN,\omega)=(\R_+\times M,d(t\theta))$ where $t>0$ is the coordinate on 
$\R_+$. Consider 
a one parameter group of symplectomorphisms $F_s:\cN\ra\cN$ of the form:
\begin{equation}
F_s(t,m)=(t,G_s(t,m)).
\label{symp}
\end{equation}
As shown in \cite{Bates} the infinitesimal generator of $F_s$ at $t=1,y=0$ is an element 
of the algebra $sp(2n+1,\R)$ of $Sp(2n+1,\R)$. 
If $X$ is vector field on $\cN$ such that 
\begin{equation}
L_X(t)=0\mbox{ and } L_X(d(t\theta))=0
\label{symp2}
\end{equation}
then its flow is a map of the form (\ref{symp}). Since 
\[L_X(d(t\theta))=d(i_X(d(t\theta)))=0\]
we deduce that (\ref{symp2}) implies:
\begin{equation}
L_X(t)=0\mbox{ and } i_Xd(t\theta)=d\widehat{h}\mbox{ for some $\widehat h:\cN\ra \R$}
\label{symp3}
\end{equation}
To see what this means, let us use Darboux coordinates. If we denote
\[X_{\widehat{h}}=\f{dt}{ds}\f{\pa}{\pa t}+
\f{dz}{ds}\f{\pa}{\pa z}+\f{dx^i}{ds}\f{\pa}{\pa x^i}+\f{dy_i}{ds}\f{\pa}{\pa y_i}\]
then (\ref{symp3}) gives:
\begin{equation}
\left\{\begin{array}{llll}
\f{dt}{ds}=0\\
\f{dx^i}{ds}=\f{1}{t}\f{\pa\widehat h}{\pa y_i}\\
\f{dy_i}{ds}=-\f{1}{t}\f{\pa\widehat h}{\pa x^i}\\
\f{dz}{ds}=-\f{\pa\widehat h}{\pa t}+\f{y_i}{t}\f{\pa\widehat h}{\pa y_i}
\end{array}
\right.
\label{Lvector}
\end{equation}
with the condition:
\begin{equation}
\f{\pa\widehat h}{\pa z}=0
\label{condZ}
\end{equation}
This has the following consequences:
\brem
Since $dt/ds=0$, the flow preserves $t=const$. Moreover
\[L_{X_{\widehat{h}}}(\widehat{h})=0,\] 
 the flow in $\cM$ associated to the vector field $X_{\widehat h}$ preserves the 
(constant) Hamiltonian function $\widehat h$.
\erem
\brem
For $t=1$ and $\widehat{h}=t\cdot h(x,y)$ one recovers the contact vector field defined in \cite{Bates}, 
that we call $X^1_h$ which is equal to:
\begin{equation}
\left\{\begin{array}{lll}
\f{dx^i}{ds}=\f{\pa h}{\pa y_i}\\
\f{dy_i}{ds}=-\f{\pa h}{\pa x^i}\\
\f{dz}{ds}=- h+y_i\f{\pa h}{\pa y_i}
\end{array}
\right.
\label{BCvector}
\end{equation}
\erem
\brem
The contactomorphisms generated by vector fields $X^1_h$ form a subgroup of the group of 
contact transformation. In fact, in \cite{Arnol'd2} V.Arnol'd defined a contact Hamiltonian
 vector field without requiring that the associated one parameter transformation 
has the form (\ref{symp}). Arnol'd's definition leads to a vector field that generates a 
contact transformation but does not preserve the contact Hamiltonian $\widehat h$.
\erem
\brem
\label{hamremark}
Consider $\widehat h_u:\cN\ra \R$ such that:
\begin{equation}
\widehat{h}_u=\f{1}{2}[-<\ga x,x>+2t<\alpha x,y>+t^2<\beta x,y>]+<v_1,y>-<v_2,x>\f{1}{t}-kt
\end{equation}
where $\alpha,\beta,\ga$ are $n\times n$ real matrices with $\beta^T=\beta$ and $\ga^T=\ga$.
One verifies that the vector field $X_{\widehat h_u}$ at $t=1$
 is the infinitesimal generator of the one parameter group in $\R^{2n+2}$
$s\ra\exp su$ with $u\in sp(2n+1,\R)$, that is, 
\begin{equation}
\left.X_{\widehat h_u}\right|_{t=1}=\left.\f{d}{ds}\right|_{s=0,t=1}\exp(su)\cdot(t,x,y,z).
\label{gen}
\end{equation}
Here
\begin{equation}
u=\left(\begin{array}{lll}
0 & 0 & 0 \\
v & A & 0 \\
k & Jv & 0
\end{array}
\right),
\end{equation}
where $J$ is the standard symplectic structure 
in $\R^{2n}$, $v\in\R^{2n}$, $k\in \R$ and $A\in sp(2n,\R)$ such that
\begin{equation}
JA=\left(\begin{array}{ll}
-\ga & \alpha  \\
\alpha^T & \beta
\end{array}
\right),
\end{equation}
Putting $(x,y)=w$, a  simple computation shows that on level set $t=1$:
\begin{equation}
\left.X_{\widehat h_u}\right|_{t=1}=(v+Aw)^T\f{\pa}{\pa w}+
(k+(Jv)^Tw)\f{\pa}{\pa z}
\label{gen2}
\end{equation}
\erem

\section{Theory of characteristics for first order PDE and  odd-symplectic vector fields}
In this section we shall follow the presentation of first order PDE given in \cite{Arnol'd}. 
Let us consider the case in which the manifold $\cM$ is the first jet bundle associated to 
$\R^n$. Locally such space 
can be modeled by a $x\in R^n$ and a germ of a function in $x$ considered up to its gradient. To each 
$m\in J^1(\R^n)$ we associate $(x,z(x),\nabla z)$ with $z(.)$ a germ of a smooth function in $\R^n$.\\
 The manifold $J^1(R^n)$ is a contact manifold with contact 1-form
\[\theta=dz-y_idx^i\]
which vanishes on any germ $z(.)$ such that $\pa z/\pa x^i=y_i$. The first jet space 
is the natural place to study the geometry of first order 
 PDEs. In fact any PDE can be understood as a submanifold of $J^1(R^n)$ (see \cite{Arnol'd}).
Let us suppose that we have a smooth function on $J^1(R^n)$:
\[h:J^1(R^n)\ra\R\mbox{ with } \f{\pa h}{\pa z}=0.\]
$h=0$ defines a manifold in $J^1(R^n)$ which corresponds to the following PDE:
\begin{equation}
h\left(x,z(x),\nabla z(x)\right)=0
\label{PDE}
\end{equation}
To equation (\ref{PDE}) we can associate a vector field $X_c$ whose integral curves satisfy 
(\ref{BCvector}). We pose the natural question 
whether $X_c$ is related to the {\it characteristic} vector field, which converts a first order PDE into  
a system of ordinary differential equations. In order to answer this question we need to 
recall the theory of characteristics. See \cite{Arnol'd} for more details.

\subsection{Theory of characteristics for first order PDE}
Let us assume that we have a first order PDE defined by the level set of 
$h:J^1(R^n)\ra \R$. In the first jet space $J^1(\R^n)$
 we can introduce {\it contact plane} 
$\Pi_x$
at any $x\in \R^n$  spanned by the vectors $V$
\begin{equation}
i_V\theta=0
\label{VV}
\end{equation}
Equation (\ref{VV}) determines a  codimension $1$ 
submanifold $E^{2n}$. Its tangent space $TE^{2n}$ can be described by the vectors $V$ such that:
\begin{equation}
i_Vdh=0
\end{equation}
\bdefi
The submanifold $E^{2n}$ is called non-characteristic, if $TE^{2n}$ is transversal to the contact plane $\Pi$. 
\edefi
\bdefi
The intersection of vector spaces
\[P_x=T_xE^{2n}\cap\Pi_x\]
 is called characteristic plane at $x$.
Therefore 
 $P$ is defined by all vectors $V$ such that
\begin{equation}
i_V\theta=0\mbox{ and }i_Vdh=0.
\label{defP}
\end{equation}
\edefi
We shall consider only the case in which $E^{2n}$ is not characteristic. Note that since $TE^{2n}$ 
and  $\Pi$ lie in $J^1(R^n)$ we have:
\[\dim P_x=\dim TE_x^{2n}+\dim\Pi_x-\dim TJ^1(R^n)=2n-1\]
The $2$-form $d\theta$, associated to the contact form $\theta$, defines a skew product in  $TJ^1(R^n)$. 
This allows us to select a vector in each 
characteristic plane $P_x$ which is skew-orthogonal to $P$. 
The characteristic vector field for (\ref{PDE}) is given by $X_c$:
\begin{enumerate}
\item $i_{X_c}\theta=0$, $X_c$ is a contact vector field, 
\item $d\theta(X_c,V)=0$ for all $V\in P$.
\label{Xcharcter}
\end{enumerate}
In symplectic geometry one can show that the skew-orthogonal complement to a $k$ dimensional vector space is 
a $2n-k$ vector space. The characteristic plane $P$ has dimension $2n-1$. Therefore 
the characteristic direction $X_c$ is uniquely determined.\\
For a given $h:J^1(R^n)\ra\R$, in local coordinates the vector field $X_c$ has components:
\begin{equation}
\left\{\begin{array}{lll}
\f{dx^i}{ds}=\f{\pa h}{\pa y_i}\\
\f{dy_i}{ds}=-\f{\pa h}{\pa x^i}-y_i\f{\pa h}{\pa z}\\
\f{dz}{ds}=y_i\f{\pa h}{\pa y_i}
\end{array}
\right.
\label{Xvector}
\end{equation}

\subsection{Cauchy Problem}
The theory of characteristics is used to solve the Cauchy problem. For an equation
\[h(x,z(x),\nabla z(x))=0\]
Cauchy data are: a manifold $\ga\subset \R^n$ and the values $z_{|_{\ga}}=\phi(x)$. We have seen that 
the equation $h=0$ corresponds 
to submanifold $E^{2n}$ in $ J^1(\R^n)$.
 Cauchy data define a submanifold $N\subset E^{n}$ called the 
{\it manifold of initial data}. It turns out that the Cauchy problem is solvable (locally) if $N$ is not characteristic,  
that is if the projection of the characteristic direction to $\R^n$
 is transversal to the tangent plane to $N$ (see 
\cite{Arnol'd}). 
Given a parametrization $x_\ga(\lambda)$ of $\ga$, it is possible to construct the family of initial data:
\[x_0=x_\ga(\lambda),z_0=\phi(x_\ga(\lambda)),p_0=p_0(\lambda)\]
where $p_0(\lambda)$ is such that:
\[h(x_\ga(\lambda),\phi(x_\ga(\lambda)),p_0(\lambda))=0\]
which forms the family of initial data for the characteristic vector field (\ref{Xvector}).
\brem
The characteristic 
vector field (\ref{Xvector}) does not depend on whether the PDE is given by $h=0$ or $h=c$ with $c\neq 0$, but 
the Cauchy problem does. In fact the value of the level set of $h$ is contained in $p_0$.
\erem

\section{The main result}
We can state:
\bth
\label{main}
For all the PDEs represented in $J^1(R^n)$ by a smooth function $h:J^1(R^n)\ra\R$ of the form:
\[h(x,\nabla z(x))=0,\]
the odd-symplectic Hamiltonian vector field (\ref{BCvector}) commutes with
the characteristic vector field (\ref{Xvector}). Moreover on $h=0$, the two vector fields
coincide.
\eth
\bpf
Consider the function $h:J^1(\R^n)\ra\R$ as a contact Hamiltonian. Then by inspection one finds that 
the infinitesimal generator $X_h^1$ can be written as follows:
\begin{equation}
X^1_h=X_c-h\f{\pa}{\pa z}.
\end{equation}
Since $h$ dos not depend on $z$, therefore $X_c(h)=0$. Also the  Lie bracket:
\[[X_c,X_h^1]=-\left[X_c,h\f{\pa}{\pa z}\right]=0\]
vanishes.\\
Now we want to compare $X^1_h$ (\ref{BCvector})  and $X_c$ (\ref{Xvector}). We notice that 
if $h:J^1(R^n)\ra\R$ does not depend on $z$, the first two components are equal.
The component ${dz}/{ds}$ differs because of the term $-h$ in (\ref{BCvector}).
We have 
\[L_{X^1_h}(h)=0\]
that is $h$ is preserved by the flow of $X^1_h$. Therefore in the manifold $h=0$
 the vector field $X^1_h$ is equal to the characteristic vector field $X_c$  
\epf

\subsection{First order PDE and $sp(2n+1,\R)$}
Let us consider a particular class of first order PDEs in $\R^n$, namely,
\begin{equation}
\label{pde1}
c_{ij}\f{\pa z(x)}{\pa x^i}\f{\pa z(x)}{\pa x^j}+b_{ij}\f{\pa z(x)}{\pa x^i}x^j
+a_{ij}x^i x^j+e_i\f{\pa z(x)}{\pa x^i}+f_ix^i=h_0
\end{equation}
or more compactly:
\begin{equation}
\label{pde01}
<B\nabla z(x),\nabla z(x)>+<e,\nabla z(x)>+<f,x>=h_0
\end{equation}
Then we have:
\bprop
\label{use-prop}
For a given first order PDE of the form (\ref{pde01}), if 
\[
B\doteq\f{1}{2}\left(\begin{array}{ll}
2a & b\\
b^T & 2c
\end{array}
\right)
\]
(where $.^T$ is the transpose) is symmetric,  
then we can associate $u\in sp(2n+1,\R)$
\begin{equation}
u=\left(\begin{array}{lll}
0 & 0 & 0 \\
v & A & 0 \\
h_0 & Jv & 0
\end{array}
\right)
\label{infgen}
\end{equation}
where $k=h_0$, $Jv=(e,f)$ and $-JB=A$, to the characteristic vector field of 
(\ref{pde1}).\\
If $a=b=0$ and $v=(v_1,0)$ in $B$ then the characteristic vector field $X_c$ associated 
to (\ref{pde1}) commutes $u$. 
\eprop
\bpf
The proof is essentially based on theorem \ref{main} and on remark \ref{hamremark}. According to 
theorem \ref{main}, equation (\ref{pde01}) can be interpreted as a contact Hamiltonian 
$h(x,y)$ where $y=\nabla z(x)$, and the contact vector field on $h=0$ 
is equal to the characteristic vector field. If the conditions on the matrices $a,b,c$ are satisfied 
 the PDE  (\ref{pde01}), on $J^1(\R^n)$ gives rise to a contact Hamiltonian  of the form:
\[
\widehat{h}(x,y,1)=<w,JAw>+<Jv,w>+h_0,\mbox{ where $w=(x,y)$}
\]
According to remark \ref{hamremark}, we can associate to such Hamiltonian an 
element of $sp(2n+1,\R)$ of the form (\ref{infgen}).\\
Let us consider the Hamiltonian
\begin{equation}
\widehat{h}_u=<a x,x>+2t<bx,y>+t^2<c x,y>+<v_1,y>-<v_2,x>\f{1}{t}-h_0t
\end{equation} 
The corresponding infinitesimal generator $X_{\widehat h_u}^1$ is:
\begin{equation}
X_{\widehat h_u}^1=\xi+Z_{\widehat h}\f{\pa}{\pa z},
\end{equation}
where
\begin{equation}
Z_{\widehat h}=<v_1,y>-<v_2,x>+h_0
\end{equation}
The characteristic vector field $X_c$ of (\ref{pde01}) is 
\begin{equation}
X_{c}=\xi+Z_{c}\f{\pa}{\pa z},
\end{equation}
where
\begin{eqnarray}
\xi=2(bx+cy)\cdot\f{\pa}{\pa x}-2(ax+by)\cdot\f{\pa}{\pa y}\\
Z_c=2<bx,y>+2<cy,y>+<v_1,y>
\end{eqnarray}
Since $h$ does not depend on $z$, the commutator is given by:
\[[X_{\widehat h}^1,X_c]=(\xi(Z_c)-\xi(Z_{\widehat h}))\f{\pa}{\pa z}\]
After some algebra one finds:
\begin{eqnarray}
\xi(Z_c)-\xi(Z_{\widehat h})=4<bx,b^Ty>+2<cy,b^Ty>+<bx,v_2>+<cy,v_2>+\nonumber\\
-2<bx,ax>-2<bx,by>-4<ax,cy>-4<by,cy>
\end{eqnarray}
Last expression is identically zero for every $(x,y)$ only for $a=b=0$ and $v_2=0$. 
\epf

\subsection{Hamilton-Jacobi theory}
In certain cases in 
Hamiltonian mechanics the construction of coordinate transformations that allows to integrate explicitly 
 the equations of motion is based on the solution of Hamilton-Jacobi equation. We  
 refer the reader to \cite{Arnol'd}. Suppose 
that we have a Hamiltonian system on $(\R^{2n},\omega)$. The Cartan 
one-form, written in local coordinates, is:
\begin{equation}
\Theta=p_idq^i-H(p,q,t)dt
\end{equation}
where $H:\R^{2n}\times\R\ra \R$ in the Hamiltonian. $\Theta$ is an integral invariant for the Hamiltonian flow. Consider  
 the extended phase space $\cM=\R^{2n}\times\R$ and two set of coordinates:
\[(q,p,t)\mbox{ and } (x,y,\tau)\]
There are two functions $K(x,y,\tau)$ and $S$ such that
\begin{equation}
p_idq^i-H(p,q,t)dt=y_idx^i-K(x,y,\tau)d\tau+dS
\label{HJ0}
\end{equation}
The function $S$ is called a ``generating function'' because it allows one 
to define a canonical transformation.\\
Equation (\ref{HJ0}) can be
 used to find the canonical transformation which realizes $(q,p,t)\ra(x,y,\tau)$. 
This transformation is usually consider in the case $t=\tau$. In this case, by  
 means of Legendre transformation one can write (\ref{HJ0}) as:
\begin{equation}
p_idq^i-H(p,q,t)dt=x^idy_i-K(x,y,t)dt+dS(q,y,t).
\label{HJ1}
\end{equation}
Then one usually chooses $K=0$ and obtains:
\begin{equation}
\left\{\begin{array}{lll}
\f{\pa S}{\pa y_i}=x_i\\
\f{\pa S}{\pa q^i}=p_i\\
H\left(q,\f{\pa S}{\pa q}\right)+\f{\pa S}{\pa t}=0
\end{array}
\right.
\label{HJ3}
\end{equation}
The third equation of (\ref{HJ3}) is the Hamilton-Jacobi equation and determines 
the function $S$ which completely describes the canonical transformation $(q,p,t)\ra (x,y,t)$. The  
Hamilton-Jacobi equation is defined in $J^1(\R^{n+1})$ and fulfills the conditions of theorem
 \ref{main}.

\subsection{Eikonal equation}
In geometric optics and WKB approximation of quantum mechanics WKB approximation
 the Eikonal equation is relevant. 
This is a first order PDE given by:
\begin{equation}
\f{1}{2}\de^{ij}\f{\pa S}{\pa x^i}\f{\pa S}{\pa x^j}-N(x)=0
\end{equation}
$S$ represents the wave front and $N$ is a  function:
\[N(x)=(n(x)/c)^2,\]
where in optics is square ratio of the refraction index $n(x)$ to the speed of light $c$; 
where in WKB theory:
\[N(x)=V(x)-E,\]
namely, potential minus the in energy.\\
When in optics we consider a medium with a piece-wise constant refraction index (multilayers),  
$N(x)$ can be easily approximated by piece-wise constant function. In quantum mechanics this happens 
when 
 the potential is piece-wise constant.
In these cases proposition \ref{use-prop} can be applied in such layer. For instance if 
\[N(x)=n_k^2/c^2\mbox{ for $x$ in the $k$-layer}\]
then in $k^{th}$-layer we have  symmetry generated by
\begin{equation}
u=\left(\begin{array}{lll}
0 & 0 & 0 \\
0 & -JA & 0 \\
n_k^2/c^2 & 0 & 0
\end{array}
\right)
\end{equation}
with
\[G=\f{1}{2}\left(\begin{array}{ll}
0 & 0 \\
0 & {\bf 1}
\end{array}
\right).
\] 
A similar analysis can be carried out in the case of semi-classical description of system which have 
piece-wise constant potential e.g. multi-well potential.

\section*{Acknowledgments}
I really wish to thank Richard Cushman for his 
 useful suggestions and long discussions we had at the  
Warwick University Mathematics Institute and his careful reading 
my manuscript. I want also to thank J.Sniatycki, and L.Bates for the interesting and stimulating 
discusions. \\ 
This paper has been supported by 
Marie Curie Fellowship (Contract n. HPMFCT-2000-00541) and by 
M.A.S.I.E. European network.

\end{document}